\documentclass[reqno]{amsart}

\usepackage{hyperref}
\usepackage{amssymb}
\usepackage{dsfont}

\def\grint{\mathrm{G}}

\def\R{\mathbb{R}}

\def\N{\mathbb{N}}
\def\Z{\mathbb{Z}}
\def\T{\mathbb{T}}

\def\Id{\mathrm{Id}}

\def\hau{\mathcal{H}}
\def\leb{\mathcal{L}}

\newcommand{\diam}[1]{\ensuremath{| #1 |}}

\def\dist{\mathrm{d}}

\def\as{\mbox{a.s.}}

\theoremstyle{plain}
\newtheorem{thm}{Theorem}
\newtheorem{prp}{Proposition}
\newtheorem{lem}{Lemma}
\newtheorem{cor}{Corollary}

\theoremstyle{definition}

\theoremstyle{remark}

\title{On randomly placed arcs on the circle}
\author{Arnaud Durand}
\address{California Institute of Technology, 1200 E. California Blvd. -- MC 217-50, Pasadena, CA 91125, USA.}
\email{durand@acm.caltech.edu}


\begin{document}

\begin{abstract}
We completely describe in terms of Hausdorff measures the size of the set of points of the circle that are covered infinitely often by a sequence of random arcs with given lengths. We also show that this set is a set with large intersection.
\end{abstract}

\maketitle

\section{Introduction and statement of the results}

Let us consider a nonincreasing sequence $\ell=(\ell_{n})_{n\geq 1}$ of positive reals converging to zero. In 1956, A.~Dvoretzky~\cite{Dvoretzky:1956bs} raised the question to find a necessary and sufficient condition on the sequence $\ell$ to ensure that the whole circle $\T=\R/\Z$ is covered almost surely by a sequence of random arcs with lengths $\ell_{1}$, $\ell_{2}$, etc. To be specific, let $A(x,l)$ denote the open arc with center $x\in\T$ and length $l>0$, that is, the set of all $y\in\T$ such that $\dist(y,x)<l/2$, where $\dist$ denotes the usual quotient distance on $\T$. Then, let $(X_{n})_{n\geq 1}$ be a sequence of random points independently and uniformly distributed on $\T$ and let
\[
E_{\ell}=\limsup_{n\to\infty} A(X_{n},\ell_{n}).
\]
Dvoretzky's problem amounts to finding a necessary and sufficient condition to ensure that
\begin{equation}\label{ascover}
\as \qquad E_{\ell}=\T.
\end{equation}
This longstanding problem, along with several of its extensions, has raised the interest of, notably, P.~Billard, P.~Erd\H{o}s, J.-P.~Kahane, P.~L\'evy and B.~Mandelbrot, and has found various applications, such as the study of multiplicative processes and that of some random series of functions, see~\cite{Kahane:1985gc,Kahane:2000ao}.

In 1972, L.~Shepp completely solved Dvoretzky's problem by showing that~(\ref{ascover}) holds if, and only if,
\begin{equation}\label{seriesshepp}
\sum_{n=1}^\infty \frac{1}{n^2}\exp(\ell_{1}+\ldots+\ell_{n})=\infty.
\end{equation}
Still, several related questions remained open. For example, when the series in~(\ref{seriesshepp}) converges, it is natural to ask which proportion of the circle is actually covered by the random arcs. In other words, what is the size of the set $E_{\ell}$? A first answer may be given by computing the value of its Lebesgue measure $\leb(E_{\ell})$. This is in fact trivial, since Fubini's theorem and the Borel-Cantelli lemma directly imply that
\begin{equation}\label{valuelebesgue}
\as \qquad \leb(E_{\ell})=\begin{cases}
0 & \mbox{if } \sum_{n}\ell_{n}<\infty \\[2mm]
1 & \mbox{if } \sum_{n}\ell_{n}=\infty.
\end{cases}
\end{equation}
A typical way of refining the description of the size of the set $E_{\ell}$ is then to compute the value of its Hausdorff dimension. This has recently been done by A.-H.~Fan and J.~Wu~\cite{Fan:2004it} in the particular case where $\ell_{n}=a/n^\alpha$, for some $a>0$ and $\alpha>1$. Their result states that
\[
\forall a>0 \quad \forall\alpha>1 \quad \as \qquad \dim E_{(a/n^\alpha)}=\frac{1}{\alpha}.
\]
Corollary~\ref{cordim} below ensures that actually, for {\em any} nonincreasing sequence $\ell=(\ell_{n})_{n\geq 1}$ converging to zero, the Hausdorff dimension of the set $E_{\ell}$ is almost surely equal to
\begin{equation}\label{defsell}
s_{\ell}=\sup\bigl\{ s\in (0,1) \:\bigl|\: \sum\nolimits_{n} {\ell_{n}}^s=\infty \bigr\}=\inf\bigl\{ s\in (0,1) \:\bigl|\: \sum\nolimits_{n} {\ell_{n}}^s<\infty \bigr\},
\end{equation}
with the convention that $\sup\emptyset=0$ and $\inf\emptyset=1$. Corollary~\ref{cordim} follows from Theorem~\ref{thmsize} below, which in fact gives the value of the Hausdorff $g$-measure of the set $E_{\ell}$ for {\em any} gauge function $g$ and not only the monomial functions used to define the Hausdorff dimension, see Subsection~\ref{sizeprop}. Therefore, this theorem provides a complete description of the size of the set $E_{\ell}$ in terms of Hausdorff measures, for any sequence $\ell$.

We also show in this note that the set $E_{\ell}$ enjoys a remarkable property originally introduced by K.~Falconer~\cite{Falconer:1994hx}, namely, it is a set with large intersection. Roughly speaking, this means that $E_{\ell}$ is ``large and omnipresent'' in the circle in some strong measure theoretic sense, see Subsection~\ref{intprop}. Sets with large intersection have been shown to arise on other occasions in probability theory -- more precisely in multifractal analysis of stochastic processes, see~\cite{Durand:2007kx,Durand:2007fk} -- as well as in other fields of mathematics, namely, number theory and dynamical systems, see~\cite{Durand:2008jk,Durand:2007uq,Durand:2006uq} and the references therein.

\subsection{Size properties of the set $E_{\ell}$}\label{sizeprop}

A typical way of completely describing the size of a subset of the circle is to determine the value of its Hausdorff $g$-measure for any gauge function $g$, see~\cite{Falconer:2003oj,Rogers:1970wb}. We call a gauge function any function $g$ defined on $[0,\infty)$ that is nondecreasing near zero, enjoys $\lim_{0^+} g=g(0)=0$ and is such that $r\mapsto g(r)/r$ is nonincreasing and positive near zero. For any gauge function $g$, the Hausdorff $g$-measure of a set $F\subseteq\T$ is then defined by
\[
\hau^g(F)=\lim_{\delta\downarrow 0}\uparrow\hau^g_{\delta}(F) \qquad\mbox{with}\qquad \hau^g_{\delta}(F)=\inf_{F\subseteq\bigcup_{p} U_{p}\atop\diam{U_{p}}<\delta} \sum_{p=1}^\infty g(\diam{U_{p}}).
\]
The infimum is taken over all sequences $(U_{p})_{p\geq 1}$ of sets with $F\subseteq\bigcup_{p} U_{p}$ and $\diam{U_{p}}<\delta$ for all $p\geq 1$, where $\diam{\cdot}$ denotes diameter. Note that if $g(r)/r$ goes to infinity at zero, every nonempty open subset of $\T$ has infinite Hausdorff $g$-measure and that, otherwise, $\hau^g$ coincides up to a multiplicative constant with the Lebesgue measure on the Borel subsets of $\T$.

The size properties of $E_{\ell}$ are then completely described by the following result.

\begin{thm}\label{thmsize}
Let $\ell=(\ell_{n})_{n\geq 1}$ be a nonincreasing sequence of positive reals converging to zero and let $g$ be a gauge function. Then, with probability one, for any open subset $V$ of $\T$,
\[
\hau^g(E_{\ell}\cap V)=\begin{cases}
\hau^g(V) &\text{if } \sum_{n}g(\ell_{n})=\infty \\[1mm]
0 &\text{if } \sum_{n}g(\ell_{n})<\infty.
\end{cases}
\]
\end{thm}

Recall that the Hausdorff dimension of a nonempty set $F\subseteq\T$ is defined with the help of the monomial functions $\Id^s$ by
\[
\dim F=\sup\{s\in (0,1) \:|\: \hau^{\Id^s}(F)=\infty\}=\inf\{s\in (0,1) \:|\: \hau^{\Id^s}(F)=0\},
\]
with the same convention regarding the infimum and the supremum of the empty set as in~(\ref{defsell}). Also, it is customary to let $\dim\emptyset=-\infty$. Using Theorem~\ref{thmsize}, it is then possible to determine the value of the Hausdorff dimension of the set $E_{\ell}$, thereby generalizing the result of~\cite{Fan:2004it}.

\begin{cor}\label{cordim}
For every nonincreasing sequence $\ell=(\ell_{n})_{n\geq 1}$ of positive reals converging to zero, with probability one, $\dim E_{\ell}=s_{\ell}$, where $s_{\ell}$ is defined by~(\ref{defsell}).
\end{cor}

Theorem~\ref{thmsize} and Corollary~\ref{cordim} are proven in Sections~\ref{proofthm} and~\ref{proofcor}, respectively.

\subsection{Large intersection properties of the set $E_{\ell}$}\label{intprop}

Rigorously, the fact that $E_{\ell}$ is a set with large intersection means that it belongs to some classes $\grint^g(V)$ of subsets of the circle that we defined in~\cite[Section~5]{Durand:2007fk}. We refer to that paper, and also to~\cite{Durand:2007uq}, for a precise definition of those classes and we content ourselves with stressing on the fact that, for any gauge function $g$ and any nonempty open subset $V$ of the circle, one may define a class $\grint^g(V)$ of {\em sets with large intersection in $V$ with respect to $g$} which, among other properties, enjoys the following.

\begin{prp}\label{grintstable}
For any gauge function $g$ and any nonempty open $V\subseteq\T$,
\renewcommand{\theenumi}{\alph{enumi}}
\begin{enumerate}
\item\label{intersect} the class $\grint^g(V)$ is closed under countable intersections;
\item\label{relsizelargeint} every set $F\in\grint^g(V)$ enjoys $\hau^{\overline{g}}(F)=\infty$ for any gauge $\overline{g}$ with $\overline{g}\prec g$;
\item $\grint^g(V)=\bigcap_{\overline{g}}\grint^{\overline{g}}(V)$ where $\overline{g}$ is a gauge function enjoying $\overline{g}\prec g$;
\item\label{openincl} $\grint^g(V)=\bigcap_{U}\grint^g(U)$ where $U$ is a nonempty open subset of $V$;
\item every $G_{\delta}$-set with full Lebesgue measure in $V$ belongs to the class $\grint^g(V)$.
\end{enumerate}
\end{prp}

The notation $\overline{g}\prec g$ means that $\overline{g}/g$ monotonically tends to infinity at zero. In view of Proposition~\ref{grintstable}, every set in the class $\grint^g(V)$ has infinite Hausdorff $\overline{g}$-measure in every nonempty open subset of $V$ for any gauge function $\overline{g}\prec g$, and any countable intersection of such sets enjoys the same property. Therefore, the classes $\grint^g(V)$ provide a rigorous way of stating that a set is large and omnipresent in $V$ in a strong measure theoretic sense.

In order to describe the large intersection properties of the set $E_{\ell}$, we shall make use of the following result, which gives a simple sufficient condition for a limsup of arcs to be a set with large intersection in the circle. It may be seen as the analog for the periodic setting of the {\em ubiquity} result established in~\cite{Durand:2007uq}.

\begin{prp}\label{ubiquity}
Let $(y_{n})_{n\geq 1}$ be a sequence in $\T$ and let $(r_{n})_{n\geq 1}$ be a sequence of positive reals converging to zero. Then, for any gauge function $g$,
\[
\leb\left(\limsup_{n\to\infty} A(y_{n},2g(r_{n}))\right)=1 \qquad\Longrightarrow\qquad \limsup_{n\to\infty} A(y_{n},2r_{n})\in\grint^g(\T).
\]
\end{prp}

Proposition~\ref{ubiquity} may be interpreted as follows. Given that any gauge function $g$ is bounded below by the identity function (up to a multiplicative constant), the limsup of the arcs $A(y_{n},2r_{n})$ may be seen as a ``reduced'' version of the limsup of the arcs $A(y_{n},2g(r_{n}))$. If the latter limsup is large and omnipresent enough to contain Lebesgue-almost every point of the circle, then its reduced version is also large and omnipresent, in the weaker sense that it belongs to the class $\grint^g(\T)$. We refer to Section~\ref{proofubiquity} for a proof of Proposition~\ref{ubiquity}.

The large intersection properties of the set $E_{\ell}$ are then completely described by the following result.

\begin{thm}\label{thmint}
Let $\ell=(\ell_{n})_{n\geq 1}$ be a nonincreasing sequence of positive reals converging to zero and let $g$ be a gauge function. Then, almost surely, for any nonempty open subset $V$ of $\T$,
\[
E_{\ell}\in\grint^g(V) \qquad\Longleftrightarrow\qquad \sum\nolimits_{n} g(\ell_{n})=\infty.
\]
\end{thm}

The remainder of this note is organized as follows: Theorems~\ref{thmsize} and~\ref{thmint} are established in Section~\ref{proofthm}, Corollary~\ref{cordim} is proven in Section~\ref{proofcor}, and the proof of Proposition~\ref{ubiquity} is given in Section~\ref{proofubiquity}. Before detailing the proofs, let us mention that we shall basically only make use of the main properties of the classes $\grint^g(V)$ given by Proposition~\ref{grintstable}, the ubiquity result given by Proposition~\ref{ubiquity} and the value~(\ref{valuelebesgue}) of the Lebesgue measure of the set $E_{\ell}$. In particular, unlike the authors of~\cite{Fan:2004it}, we do not need to call upon any specific result on the spacings between the random centers $X_{n}$ of the arcs. This also means that our method can straightforwardly be extended to the case of random balls on the $d$-dimensional torus for any $d\geq 2$.

\section{Proofs of Theorems~\ref{thmsize} and~\ref{thmint}}\label{proofthm}

Theorems~\ref{thmsize} and~\ref{thmint} follow from four lemmas which we now state and prove. Throughout the section, $\ell=(\ell_{n})_{n\geq 1}$ is a nonincreasing sequence of positive reals converging to zero.

\begin{lem}\label{convhau}
For any gauge function $g$,
\[
\sum\nolimits_{n} g(\ell_{n})<\infty \qquad\Longrightarrow\qquad \forall V\mbox{ open} \quad \hau^g(E_{\ell}\cap V)=0.
\]
\end{lem}

\begin{proof}
For any $\delta>0$, there is an integer $n_{0}\geq 1$ such that $0<\ell_{n}<\delta$ for any $n\geq n_{0}$. Moreover, the set $E_{\ell}$ is covered by the arcs $A(X_{n},\ell_{n})$ for $n\geq n_{0}$, so that $\hau^g_{\delta}(E_{\ell})\leq\sum_{n=n_{0}}^\infty g(\ell_{n})$. If the series $\sum_{n} g(\ell_{n})$ converges, then letting $n_{0}$ tend to infinity and $\delta$ go to zero yields $\hau^g(E_{\ell})=0$.
\end{proof}

\begin{lem}
For any gauge function $g$,
\[
\sum\nolimits_{n} g(\ell_{n})<\infty \qquad\Longrightarrow\qquad \forall V\neq\emptyset\mbox{ open} \quad E_{\ell}\not\in\grint^g(V).
\]
\end{lem}

\begin{proof}
Let us assume that the series $\sum_{n} g(\ell_{n})$ converges. Then, one may build a gauge function $\overline{g}$ such that $\overline{g}\prec g$ and the series $\sum_{n} \overline{g}(\ell_{n})$ converges too, for example by adapting a construction given in the proof of~\cite[Theorem~3.5]{Durand:2007lr}. By Lemma~\ref{convhau}, the set $E_{\ell}$ has Hausdorff $\overline{g}$-measure zero in $V$ and thus cannot belong to the class $\grint^g(V)$, due to Proposition~\ref{grintstable}(\ref{relsizelargeint}).
\end{proof}

\begin{lem}\label{divgrint}
For any gauge function $g$,
\[
\sum\nolimits_{n} g(\ell_{n})=\infty \qquad\Longrightarrow\qquad \as\quad\forall V\neq\emptyset\mbox{ open} \quad E_{\ell}\in\grint^g(V).
\]
\end{lem}

\begin{proof}
If the series $\sum_{n} g(\ell_{n})$ diverges, then $\sum_{n} g(\ell_{n}/2)$ diverges as well (because $r\mapsto g(r)/r$ is nonincreasing near zero). Hence, thanks to~(\ref{valuelebesgue}), the limsup of the arcs $A(X_{n},2g(\ell_{n}/2))$ has Lebesgue measure one with probability one. We conclude using Proposition~\ref{ubiquity} and Proposition~\ref{grintstable}(\ref{openincl}).
\end{proof}

\begin{lem}
For any gauge function $g$,
\[
\sum\nolimits_{n} g(\ell_{n})=\infty \qquad\Longrightarrow\qquad \as\quad\forall V\mbox{ open} \quad \hau^g(E_{\ell}\cap V)=\hau^g(V).
\]
\end{lem}

\begin{proof}
We may obviously assume that $V$ is nonempty. Let us suppose that the series $\sum_{n} g(\ell_{n})$ diverges. Then, again by following a construction given in the proof of~\cite[Theorem~3.5]{Durand:2007lr}, it is possible to build a gauge function $\underline{g}$ such that $g\prec\underline{g}$ and the series $\sum_{n} \underline{g}(\ell_{n})$ diverges too, provided that $g\prec\Id$. Therefore, thanks to Lemma~\ref{divgrint}, the set $E_{\ell}$ belongs to the class $\grint^{\underline{g}}(V)$. Hence, $\hau^g(E_{\ell}\cap V)=\infty=\hau^g(V)$, owing to Proposition~\ref{grintstable}. In the case where $g\not\prec\Id$, the Hausdorff $g$-measure coincides, up to a multiplicative constant, with the Lebesgue measure on the Borel subsets of the circle and the result follows from~(\ref{valuelebesgue}).
\end{proof}

\section{Proof of Corollary~\ref{cordim}}\label{proofcor}

In order to prove Corollary~\ref{cordim}, let us consider a nonincreasing sequence $\ell=(\ell_{n})_{n\geq 1}$ of positive reals converging to zero. Theorem~\ref{thmsize}, along with the definition~(\ref{defsell}) of the real $s_{\ell}$, implies that for any real $s\in (0,1)$, with probability one,
\[
\hau^{\Id^s}(E_{\ell})=\begin{cases}
\infty &\text{if } s<s_{\ell} \\[1mm]
0 &\text{if } s>s_{\ell}.
\end{cases}
\]
Let us assume that $s_{\ell}\in (0,1]$. Then, for all $m$ large enough, with probability one, the set $E_{\ell}$ has infinite Hausdorff $\Id^{s_{\ell}-1/m}$-measure, so that its Hausdorff dimension is at least $s_{\ell}-1/m$. Therefore, the dimension of $E_{\ell}$ is almost surely at least $s_{\ell}$. Likewise, if $s_{\ell}\in [0,1)$, then $E_{\ell}$ has $\Id^{s_{\ell}+1/m}$-measure zero with probability one for all $m$ large enough, so that its Hausdorff dimension is almost surely at most $s_{\ell}$. As a result, with probability one, $\dim E_{\ell}=s_{\ell}$ if $s_{\ell}\in (0,1]$ and $\dim E_{\ell}\leq 0$ if $s_{\ell}=0$.

It remains to establish that $E_{\ell}$ is almost surely nonempty when $s_{\ell}=0$. Note that the set $E_{(1/n)}$, obtained by picking $\ell_{n}=1/n$, has Lebesgue measure one with probability one, by virtue of~(\ref{valuelebesgue}). Let us assume that this property holds. Furthermore, note that $\ell_{n}={\rm o}(1/n)$ as $n$ goes to infinity, thanks to Olivier's theorem~\cite{Olivier:1827eu}. In particular, $\ell_{n}\leq 1/n$ for any integer $n$ greater than or equal to some $n_{1}\geq 1$. Let $I_{1}=A(X_{n_{1}},\ell_{n_{1}}/2)$. The union over $n>\max\{n_{1},8/\ell_{n_{1}}\}$ of the arcs $A(X_{n},1/n)$ has full Lebesgue measure in the circle, so its intersection with the arc $A(X_{n_{1}},\ell_{n_{1}}/4)$ is nonempty. Therefore, there is an integer $n_{2}>\max\{n_{1},8/\ell_{n_{1}}\}$ such that $A(X_{n_{2}},1/n_{2})\subseteq I_{1}$. Then, let $I_{2}=A(X_{n_{2}},\ell_{n_{2}}/2)$. Repeating this procedure, one may obtain a nested sequence of open arcs $I_{n}$ and the intersection of their closures yields a point that belongs to the set $E_{\ell}$.

\section{Proof of Proposition~\ref{ubiquity}}\label{proofubiquity}

Let $g$ be a gauge function such that the limsup of the arcs $A(y_{n},2g(r_{n}))$ has Lebesgue measure one. Thus, following the terminology of~\cite{Durand:2007uq}, the family $(k+\dot y_{n},g(r_{n}))_{(k,n)\in\Z\times\N}$ is a homogeneous ubiquitous system in $\R$. Here, each $\dot y_{n}$ is the only real in $[0,1)$ such that $\phi(\dot y_{n})=y_{n}$, where $\phi$ denotes the canonical surjection from $\R$ onto $\T$. Thanks to~\cite[Theorem~2]{Durand:2007uq}, the set of all reals $x$ such that $|x-k-\dot y_{n}|<r_{n}$ for infinitely many $(k,n)\in\Z\times\N$ belongs to the class $\grint^g(\R)$ of sets with large intersection in $\R$ with respect to the gauge function $g$, which is defined in~\cite{Durand:2007uq}. Equivalently, the inverse image under $\phi$ of the limsup of the arcs $A(y_{n},2r_{n})$ belongs to $\grint^g(\R)$, which ensures that this limsup belongs to the class $\grint^g(\T)$, see~\cite[Section~5]{Durand:2007fk}.

\end{document}